\documentclass[graybox]{svmult}

\usepackage{type1cm}        
\usepackage{makeidx}         
\usepackage{graphicx}        
\usepackage{multicol}        
\usepackage[bottom]{footmisc}

\usepackage{amsmath}
\usepackage{amsbsy}
\usepackage{amssymb}
\usepackage{amsfonts}
\usepackage{mathrsfs}
\usepackage{graphicx,color}
\usepackage{bm}
\usepackage{bbm}
\usepackage{booktabs}
\usepackage{here} 
\usepackage{enumerate}
\usepackage{newtxtext}       %
\usepackage{newtxmath}       
\usepackage{microtype}

\allowdisplaybreaks

\makeindex             

\newcommand{\weakcondWinP}{\,\overset{\mathrm{P}}{\underset{W}{\rightsquigarrow}}\,}
\newcommand{\weakcondVinP}{\,\overset{\mathrm{P}}{\underset{V}{\rightsquigarrow}}\,}

\newcommand{\weakcondxiinP}{\,\overset{\mathrm{P}}{\underset{\xi}{\rightsquigarrow}}\,}

\newcommand{\inP}{\overset{\mathrm{P}}{\longrightarrow}}

\newcommand{\tendsto}{\to}
\newcommand{\indic}{\mathbbm{1}}
\renewcommand{\qed}{\hfill$\Box$}
\renewcommand{\Pr}{\operatorname{P}}

\renewcommand{\E}{\mathrm{E}}
\newcommand{\Var}{\operatorname{Var}}
\newcommand{\Cov}{\operatorname{Cov}}
\renewcommand{\d}{\mathrm{d}}

\newcommand{\vecR}{\boldsymbol{R}}

\newcommand{\vecV}{\boldsymbol{V}}
\newcommand{\vecX}{\boldsymbol{X}}

\newcommand{\vecu}{\boldsymbol{u}}
\newcommand{\vecv}{\boldsymbol{v}}

\newcommand{\vecx}{\boldsymbol{x}}

\newcommand{\vc}[1]{\bm{#1}}
\newcommand{\half}{\frac{1}{2}}

\renewcommand{\epsilon}{\varepsilon}
\newcommand{\eps}{\varepsilon}

\newcommand{\bbC}{\mathbb{C}}

\newcommand{\bbF}{\mathbb{F}}
\newcommand{\bbG}{\mathbb{G}}

\newcommand{\bbR}{\mathbb{R}}

\newcommand{\bbU}{\mathbb{U}}

\newcommand{\tilbbC}{\tilde{\mathbb{C}}}

\renewcommand{\I}{\mathbbm{1}}

\renewcommand{\ge}{\geqslant}
\renewcommand{\geq}{\ge}
\renewcommand{\le}{\leqslant}
\renewcommand{\leq}{\le}

\newcommand{\diff}{\mathrm{d}}

\newcommand{\cop}{\mathbb{C}_n}

\newcommand{\copb}{\mathbb{C}_n^{\beta}}
\newcommand{\copbst}{\mathbb{C}_n^{\beta *}}
\newcommand{\copbh}{\mathbb{C}_n^{\beta \#}}

\theoremstyle{plain}
\theorembodyfont{\upshape}      
\newtheorem{Def}{Definition}[section]

\newtheorem{Cond}[Def]{Condition}
\newtheorem{Alg}[Def]{Algorithm}

\theoremstyle{break}

\theoremstyle{plain}
\theorembodyfont{\slshape}

\newtheorem{Prop}[Def]{Proposition}

\newtheorem{Thm}[Def]{Theorem}

\theoremstyle{break}

\definecolor{auburn}{rgb}{0.43, 0.21, 0.1}
\definecolor{britishracinggreen}{rgb}{0.0, 0.26, 0.15}
\definecolor{burntumber}{rgb}{0.54, 0.2, 0.14}
\definecolor{carmine}{rgb}{0.59, 0.0, 0.09}
\definecolor{aurometalsaurus}{rgb}{0.43, 0.5, 0.5}
\definecolor{darkcandyapplered}{rgb}{0.64, 0.0, 0.0}
\definecolor{darkpowderblue}{rgb}{0.0, 0.2, 0.6}
\definecolor{darkraspberry}{rgb}{0.53, 0.15, 0.34}
\definecolor{gray}{rgb}{0.4, 0.4, 0.4 }

\begin{document}

\title*{Resampling Procedures with Empirical Beta Copulas}
\author{Anna Kiriliouk, Johan Segers, and Hideatsu Tsukahara}
\institute{Anna Kiriliouk \at Facult\'{e} des sciences \'{e}conomiques, sociales et de gestion, 
Universit\'{e} de Namur, Rue de Bruxelles 61, B-5000 Namur, Belgium, \email{anna.kiriliouk@unamur.be}
\and Johan Segers \at Institut de Statistique, Biostatistique et Sciences Actuarielles, 
Universit\'{e} catholique de Louvain, Voie du Roman Pays 20, B-1348 Louvain-la-Neuve, Belgium, 
\email{johan.segers@uclouvain.be}
\and Hideatsu Tsukahara \at Faculty of Economics, Seijo University, 6--1--20 Seijo, Setagaya-ku, 
Tokyo, 157-8511, Japan, \email{tsukahar@seijo.ac.jp}}
%
%
\maketitle


\abstract{%
	The empirical beta copula is a simple but effective smoother of the empirical copula. Because it is a genuine copula, from which it is particularly easy to sample, it is reasonable to expect that resampling procedures based on the empirical beta copula are expedient and accurate. In this paper, after reviewing the literature on some bootstrap approximations for the empirical copula process, we first show the asymptotic equivalence of several bootstrapped processes related to the empirical and empirical beta copulas.  Then we investigate the finite-sample properties of resampling schemes based on the empirical (beta) copula by Monte Carlo simulation. More specifically, we consider interval estimation for functionals such as the rank correlation coefficients and dependence parameters of several well-known families of copulas. Here we construct confidence intervals using several methods, and compare their accuracy and efficiency. We also compute the actual size and power of symmetry tests based on several resampling schemes for the empirical and empirical beta copulas.%
}

\section{Introduction}

Let $\vecX_i=(X_{i1},\ldots,X_{id})$, $i\in\{1,\ldots,n\}$, be independent and identically 
distributed (i.i.d.) random vectors, and assume that the cumulative distribution function, $F$, of $\vecX_i$
is continuous.  By Sklar's theorem \cite{Sklar59}, there exists a unique copula, $C$, 
such that
\[
  F(\vecx)=C\bigl(F_1(x_1),\ldots,F_d(x_d)\bigr),\quad \vecx=(x_1,\ldots,x_d)\in\bbR^d, 
\]
where $F_j$ is the $j$th marginal distribution function of $F$.  In fact, in the continuous case, 
we have $C(\vecu)=F\bigl(F_1^-(u_1),\ldots,F_d^-(u_d)\bigr)$ for $\vecu=(u_1,\ldots,u_d) \in [0,1]^d$, where $H^-(u)=\inf\{t\in\bbR\colon H(t)\geq u\}$ is the generalized inverse of a distribution 
function $H$.  The empirical copula $\bbC_n$ \cite{Deheu79} is defined by
\begin{equation*}
  \bbC_n(\vecu):=\bbF_n\bigl(\bbF_{n1}^-(u_1),\ldots,\bbF_{nd}^-(u_d)\bigr),
\end{equation*}
where, for $j\in\{1,\ldots,d\}$, 
\[
  \bbF_n(\vecx):=\frac{1}{n}\sum_{i=1}^n\mathbbm{1} \left\{ X_{i1} \leq x_1,\ldots, 
    X_{id} \leq x_d \right\}, \quad
  \bbF_{nj}(x_j):=\frac{1}{n}\sum_{i=1}^n\mathbbm{1} \left\{ X_{ij} \leq x_j\right\}.
\]

For $i\in\{1,\ldots,n\}$ and $j\in\{1,\ldots,d\}$, let $R_{ij,n}$ be the rank of $X_{ij}$ 
among $X_{1j},\ldots,X_{nj}$; namely, 
\begin{equation}
\label{eq:ranks}
  R_{ij,n} = \sum_{k=1}^n \indic\{ X_{kj} \le X_{ij}\}.
\end{equation}
A frequently used rank-based version of the empirical copula is given by
\begin{equation}
\label{eq:cop:rank}
  \tilbbC_n(\vecu):=\frac{1}{n}\sum_{i=1}^n\prod_{j=1}^d\indic\biggl\{
  \frac{R_{ij,n}}{n}\leq u_j\biggr \}.
\end{equation}
In the absence of ties, we have 
\begin{equation}
\label{ineq:emp-cop:rank-cop}
  \lVert \tilbbC_n-\bbC_n \rVert_\infty
  :=\sup_{\vecu\in [0,1]^d} \lvert \tilbbC_n(\vecu)-\bbC_n(\vecu) \rvert
  \leq \frac{d}{n}. 
\end{equation}
The functions $\bbC_n$ and $\tilbbC_n$ are both piecewise constant and cannot be genuine copulas. When the sample size is small, they suffer from the presence of ties when used in resampling.

The empirical beta copula \cite{SST2017}, defined in Section~\ref{sec:beta}, is a simple but effective way of correcting 
and smoothing the empirical copula. Even though its asymptotic distribution is the same as that of the usual empirical copula, its accuracy in small samples is usually better, partly because it is itself always a genuine copula. Moreover, drawing random samples from the empirical beta copula is quite straightforward. 

Because of these properties, it is reasonable to expect that simple and accurate resampling schemes for the empirical copula process can be constructed based on the empirical beta copula.
For tail copulas, which are limit functions describing the asymptotic behavior of a copula in the corner of the unit cube, 
a simulation study in \cite{Kiriliouk-Segers-Tafakori2018} shows that a bootstrap based on the empirical beta copula performs significantly better than the direct multiplier bootstrap of \cite{Buech-Dette2010}.  
The purpose of this study is to further investigate the finite-sample and the asymptotic behavior 
of this resampling method for general copulas.

The paper is structured as follows. In Section~\ref{sec:bs-review}, we review and discuss the literature on resampling methods for the empirical copula process.  The asymptotic properties
of two resampling procedures based on the empirical beta copula are investigated in Section~\ref{sec:beta}.
In Section~\ref{sec:sim}, extensive simulation studies are conducted to demonstrate the effectiveness of resampling procedures based on the empirical beta copula to construct confidence intervals for several copula functionals, and to test the shape constraints on the copula.  We conclude the paper with some discussion and open questions in Section~\ref{sec:concl-rmks}.
All proofs are relegated to the Appendix.

\section{Review on bootstrapping empirical copula processes}
\label{sec:bs-review}

In this section, we give a short review on bootstrapping empirical copula processes, incorporating several newer improvements. We limit ourselves to i.i.d.\ sequences. Note that extensions to stationary time series are considered in \cite{Buech-Volg2013}, among others.

First, we recall a basic result on the weak convergence of the empirical copula process.  
Let $\ell^\infty([0,1]^d)$ be the Banach space of real-valued, bounded functions on $[0,1]^d$, 
equipped with the supremum norm $\lVert\,\cdot\,\rVert_\infty$.
The arrow $\rightsquigarrow$ denotes weak convergence in the sense used in \cite{Vaart-Wellner}.  
The following is the only condition needed for our convergence results.
\begin{Cond}
\label{cond:partial-deriv}
For each $j\in\{1,\ldots,d\}$, the copula $C$ has a continuous first-order
partial derivative $\dot{C}_j(\vecu)=\partial C(\vecu)/\partial u_j$ on the set 
$\{\vecu\in [0,1]^d\colon 0<u_j<1\}$. 
\end{Cond}
The next theorem is proved in \cite{Segers2012}. Let $\bbU^C$ denote a $C$-pinned Brownian sheet, 
that is, a centered Gaussian process on $[0, 1]^d$ with continuous trajectories and covariance function
\begin{equation}
\label{eq:UC}
\Cov \bigl\{ \bbU^C(\vecu),\bbU^C(\vecv) \bigr\} 
= C(\vecu\wedge\vecv)-C(\vecu) \, C(\vecv),
\qquad \vecu, \vecv \in [0, 1]^d.
\end{equation}

\begin{Thm}
\label{thm:conv-emp-cop}
Suppose Condition~\ref{cond:partial-deriv} holds.  Then we have
\[
  \bbG_n:=\sqrt{n}(\bbC_n-C)\rightsquigarrow\bbG^C, \qquad n \to \infty,
\]
in $\ell^\infty([0,1]^d)$, where
\begin{equation*}
  \bbG^C(\vecu) := \bbU^C(\vecu)-\sum_{j=1}^d \dot{C}_j(\vecu)\, \bbU^C(\mathbf{1},u_j,\mathbf{1}),
\end{equation*}
with $u_j$ appearing at the $j$th coordinate.
\end{Thm}

Next, we introduce notation for the convergence of conditional laws in probability given the data as defined in \cite{Kosorok2008}; see also \cite[Section~2.9]{Vaart-Wellner}.  
Let 
\begin{multline}
\label{eq:BL1}
  \operatorname{BL}_1:=\{h\colon\ell^\infty([0,1]^d)\tendsto \bbR
  \mid \lVert h\rVert_\infty\leq 1\ 
  \text{and}\ \lvert h(x)-h(y) \rvert \leq \lVert x-y\rVert_\infty \\
   \text{for all}\ x,y\in \ell^\infty([0,1]^d)\}.
\end{multline}
If $\hat{X}_n$ is a sequence of bootstrapped processes in $(\ell^\infty([0,1]^d), 
\lVert\,\cdot\,\rVert_\infty)$ with random weights $W$, then the notation
\begin{equation}
\label{eq:weakcondWinP}
  \hat{X}_n\weakcondWinP X, \qquad n\tendsto\infty
\end{equation}
means that 
\begin{equation}
\label{eq:weakcondWinP:2}
	\left.
	\begin{split}
  \sup_{h\in{\operatorname{BL}}_1} \lvert \E_W [h(\hat{X}_n)] - \E[h(X)] \rvert 
  &\longrightarrow 0
	\quad \text{in outer probability}, \\
  \E_W[h(\hat{X}_n)^*]-\E_W[h(\hat{X}_n)_*]
  &\inP 0 \quad\mbox{for all}\ h\in\operatorname{BL}_1. 
  \end{split}
  \right\}
\end{equation}
Here, the notation $\E_W$ indicates conditional expectation over the weights $W$ given
the data $\vecX_1, \ldots, \vecX_n$, and $h(\hat{X}_n)^*$ and $h(\hat{X}_n)_*$ denote
the minimal measurable majorant and maximal measurable minorant, respectively,
with respect to the joint data $\vecX_1, \ldots, \vecX_n, W$.

In the following, the random weights $W$ can signify different things: a multinomial random vector when drawing from the data with replacement; i.i.d.\ multipliers in the multiplier bootstrap; or vectors of order statistics from the uniform distribution when resampling from the empirical beta copula. In \eqref{eq:weakcondWinP}, the symbol $W$ will then be changed accordingly.

\subsection{Straightforward bootstrap}
\label{subsec:straight-bs}
Let $(W_{n1},\ldots,W_{nn})$ be a multinomial random vector with probabilities $(1/n,\ldots,1/n)$, 
independent of the sample $\vecX_1, \ldots, \vecX_n$.  Set
\begin{equation*}
  \bbC_n^*(\vecu)
  = 
  \bbF_n^* \bigl(\bbF_{n1}^{*-}(u_1),\ldots,\bbF_{nd}^{*-}(u_d) \bigr),
\end{equation*}
where 
\begin{align*}
  \bbF_n^*(\vecx)&:=\frac{1}{n}\sum_{i=1}^nW_{ni} \prod_{j=1}^d \indic\{X_{ij}\leq x_j\},\\
  \bbF_{nj}^*(x_j)&:=\frac{1}{n}\sum_{i=1}^nW_{ni}\indic\{X_{ij}\leq x_j\},\quad j\in\{1,\ldots,d\}.
\end{align*}
We can also define the bootstrapped version of the rank-based empirical copula 
\begin{equation}
\label{copboot2}
  \tilbbC_n^{*}(\vecu)=  
  \frac{1}{n} \sum_{i=1}^n W_{ni} \prod_{j=1}^d 
  \I \left\{ \frac{R_{ij,n}^*}{n} \leq u_j \right\},
\end{equation}
where
\begin{equation}
\label{eq:bootstrappedRanks}
  R^*_{ij,n} 
  = \sum_{k=1}^n W_{nk} \I \left\{ X_{kj} \leq X_{ij} \right\}.
\end{equation}
Since a bootstrap sample will have ties with a (large) positive probability, 
the bound \eqref{ineq:emp-cop:rank-cop} is no longer valid for $\bbC_n^{*}$ and $\tilbbC_n^{*}$.
However, we can prove the following.
\begin{Prop}
\label{prop:bs-emp-cop-diff}
\begin{equation}
\label{eq:copEquivCop}
  \lVert \bbC_n^{*} - \tilbbC_n^{*} \rVert_\infty = O_p \bigl( n^{-1} \log n \bigr),
  \qquad n \to \infty.
\end{equation}
\end{Prop}
The proof of Proposition~\ref{prop:bs-emp-cop-diff} is given in the Appendix. Convergence in probability of the conditional laws
\begin{equation*}
  \sqrt{n}(\bbC_n^*-\bbC_n) \weakcondWinP \bbG^C, \quad n\tendsto\infty,
\end{equation*}
in the space $\ell^\infty([0,1]^d)$ is shown in \cite{Fer-Rad-Weg04} under the condition that all partial derivatives $\dot{C}_j$ exist and are continuous on $[0,1]^d$, and in~\cite{Buech-Volg2013} under the weaker Condition~\ref{cond:partial-deriv}. From \eqref{ineq:emp-cop:rank-cop} and Proposition~\ref{prop:bs-emp-cop-diff}, we also have
\begin{equation}
\label{eq:alphatildeconv}
  \tilde{\alpha}_n := \sqrt{n} (\tilde{\bbC}_n^* - \tilde{\bbC}_n) \weakcondWinP \bbG^C, \quad n\tendsto\infty.
\end{equation}

\subsection{Multiplier bootstrap with estimated partial derivatives}

The multiplier bootstrap for the empirical copula, proposed by \cite{Remillard-Scaillet2009}, has proved useful in many problems. In \cite{Buech-Dette2010}, this method is found to exhibit better finite-sample performance than other resampling methods for the empirical copula process.  Here, we present a modified version proposed by \cite{Buech-Dette2010}, which we employ in the simulation studies in Section~\ref{sec:sim}.

Let $\xi_1,\ldots,\xi_n$ be i.i.d.\ nonnegative random variables, independent of the data, 
with $\E(\xi_i)=\mu$, $\Var(\xi_i)=\tau^2>0$, and $\|\xi_i\|_{2,1}:=\int_0^\infty\sqrt{\Pr(|\xi_i|>x)}\,{\d}x<\infty$. 
Put $\overline{\xi}_n:=n^{-1}\sum_{i=1}^n\xi_i$, and set
\begin{align*}
  \bbC_n^{\circ}(\vecu) &:=\frac{1}{n} \sum_{i=1}^n 
  \frac{\xi_i}{\overline{\xi}_n} \prod_{j=1}^d \I \left\{X_{ij} \leq \bbF_{nj}^-(u_j)\right\},\\
  \tilbbC_n^\circ(\vecu) &:= \frac{1}{n} \sum_{i=1}^n \frac{\xi_i}{\overline{\xi}_n} 
  \prod_{j=1}^d \I \left\{\bbF_{nj}(X_{ij}) \leq u_j \right\}.
\end{align*}
Define $\beta_n^\circ := \sqrt{n}(\mu/\tau)(\bbC_n^{\circ}-\bbC_n)$  
and $\tilde{\beta}_n^\circ := \sqrt{n}(\mu/\tau)(\tilbbC_n^{\circ}-\tilbbC_n)$.
Using Theorem~2.6 in \cite{Kosorok2008} and the almost sure convergence $\lVert \bbF_{nj}^{-}-I\rVert_\infty\tendsto 0$, 
where $I$ is the identity function on $[0,1]$, we can show that
\begin{equation*} 
	\beta_n^\circ \weakcondxiinP \bbU^C 
	\qquad\text{and}\qquad
	\tilde{\beta}_n^\circ \weakcondxiinP \bbU^C,
	\qquad
	n \to \infty.
\end{equation*}
Hence, if $\hat{\dot{C}}_j(\vecu)$ is an estimate for $\dot{C}_j(\vecu)$, where finite differencing is applied to the empirical copula at a spacing proportional to $n^{-1/2}$, then the processes
\begin{equation*}
	\left\{
	\begin{split}
  \alpha_n^{\mathrm{pdm}\circ}(\vecu)
  &:=\beta_n^\circ(\vecu)- \textstyle\sum_{j=1}^d 
  \hat{\dot{C}}_j(\vecu) \, \beta_n^\circ(\mathbf{1},u_j,\mathbf{1}) \\
  \tilde{\alpha}_n^{\mathrm{pdm}\circ}(\vecu)
  &:=\tilde{\beta}_n^\circ(\vecu)- \textstyle\sum_{j=1}^d 
  \hat{\dot{C}}_j(\vecu) \, \tilde{\beta}_n^\circ(\mathbf{1},u_j,\mathbf{1})
  \end{split}
  \right.
\end{equation*}
yield \textit{conditional} approximations of $\bbG^C$, where ``pdm'' stands for ``partial derivatives
multiplier''.  That is, we have
\[
  \alpha_n^{\mathrm{pdm}\circ}\weakcondxiinP \bbG^C \quad \mbox{and}\quad
  \tilde{\alpha}_n^{\mathrm{pdm}\circ}\weakcondxiinP \bbG^C, 
  \qquad n\tendsto\infty.
\]
%

\section{Resampling with empirical beta copulas}
\label{sec:beta}

The \emph{empirical beta copula} \cite{SST2017} is defined as
\begin{equation*}
  \bbC_n^\beta(\vecu)=\frac{1}{n}\sum_{i=1}^n\prod_{j=1}^d F_{n,R_{ij,n}}(u_j), \qquad
  \vecu \in [0,1]^d,
\end{equation*}
where $R_{ij,n}$ denotes the rank, as in \eqref{eq:ranks}, and where, for $u\in [0,1]$ and $r \in \{1, \ldots, n\}$,
\begin{equation}
\label{eq:betacdf}
  F_{n,r}(u) = \sum_{s=r}^n \binom{n}{s} u^s (1-u)^{n-s}
\end{equation}
is the cumulative distribution function of the beta distribution $\mathcal{B}(r,n+1-r)$. 
Note that $\Pr(U \leq u) = \Pr(S \geq r)$, for $U \sim \mathcal{B}(r,n+1-r)$ and $S \sim \operatorname{Bin}(n, u)$.  
In this section, we examine the asymptotic properties of two resampling procedures based on the empirical beta copula.

\subsection{Standard bootstrap for the empirical beta copula}

Let $(W_{n1},\ldots,W_{nn})$ be a multinomial random vector with success probabilities $(1/n,\ldots,1/n)$, independent of the original sample. Set
\begin{equation*}
  \copbst ( \vc{u}) 
  = 
  \frac{1}{n} \sum_{i=1}^n W_{ni} \prod_{j=1}^d F_{n,R^*_{ij,n}} (u_j),
\end{equation*}
where $R^*_{ij,n}$ are the bootstrapped ranks in \eqref{eq:bootstrappedRanks}.  
Let $S_j \sim \operatorname{Bin}(n, u_j)$, for $j = 1, \ldots, d$, be $d$ independent binomial random variables. Let $\E_S$ denote the expectation with respect to $(S_1, \ldots, S_d)$, conditional on the sample and the multinomial random vector. 
It then follows that
\begin{align*}
  \copbst ( \vc{u}) 
  &=
  \frac{1}{n} \sum_{i=1}^n W_{ni} \prod_{j=1}^d
  \E_S\left[ \I \left\{ \frac{R^*_{ij,n}}{n} \le \frac{S_j}{n}  \right\} \right] 
  =
  \E_S \left[ \tilbbC_n^{*}(S_1/n, \ldots, S_d/n) \right],
\end{align*}
where $\tilbbC_n^*$ is the bootstrapped rank-based empirical copula in \eqref{copboot2}. Similarly, the empirical beta copula is
\[
  \copb ( \vc{u}) 
  = \frac{1}{n} \sum_{i=1}^n \prod_{j=1}^d F_{n,R_{ij,n}} (u_j) 
  = \E_S \left[ \tilbbC_n( S_1/n, \ldots, S_d/n ) \right],
\]
where $\tilbbC_n$ is the rank-based empirical copula in \eqref{eq:cop:rank}. Consider the bootstrapped processes $\tilde{\alpha}_n$ defined in \eqref{eq:alphatildeconv} and $\alpha_n^\beta := \sqrt{n} ( \copbst - \bbC_n^\beta)$.
We find
\begin{equation}
\label{eq:alphaES}
  \alpha_n^\beta(\vc{u})
  =
  \E_S[ \tilde{\alpha}_n(S_1/n, \ldots, S_d/n) ].
\end{equation}
Using the weak convergence of the bootstrapped process $\tilde{\alpha}_n$, we prove the following proposition.
As a result, the consistency of the bootstrapped process $\tilde{\alpha}_n$ of the (rank-based) empirical copula in \eqref{eq:alphatildeconv} entails consistency of the one for the empirical beta copula.

\begin{Prop}
\label{prop:asymp-equiv}
Under Condition~\ref{cond:partial-deriv}, we have 
\begin{equation}
\label{eq:alphabeta:conv}
  \sup_{\vc{u} \in [0, 1]^d}
  \lvert \alpha_n^\beta(\vc{u}) - \tilde{\alpha}_n(\vc{u}) \rvert
  =
  o_p(1),
  \qquad n \to \infty,
\end{equation}
and thus $\alpha^\beta_n \weakcondWinP \bbG^C$ as $n\tendsto\infty$.
\end{Prop}

\subsection{Bootstrap by drawing samples from the empirical beta copula}

The original motivation of \cite{SST2017} was resampling; the uniform random variables generated independently and rearranged in the order specified by the componentwise ranks of the original sample might, in some sense, be considered a bootstrap sample.  Although this idea turned out to be not entirely correct, it still led to the discovery of the empirical beta copula.  In the same spirit, it is natural to study the bootstrap method using samples drawn from the empirical beta copula $\bbC_n^\beta$.

It is in fact very simple to generate a random variate $\vecV$ from $\bbC_n^\beta$. 

\begin{svgraybox}\vspace{-\baselineskip}
\begin{Alg}
\label{alg:beta-sampling}
Given the ranks $R_{ij,n}=r_{ij}, \; j=1,\ldots,d$, of the original sample: 
\begin{enumerate}[1.]
\setlength{\itemsep}{0pt}
\item Generate $I$ from the discrete uniform distribution on $\{1,\ldots, n\}$.
\item Independently generate $V_{j}^{\#}\sim \mathcal{B}(r_{Ij},n+1-r_{Ij})$,\ \,$j\in\{1,\ldots,d\}$.
\item Set $\vecV^{\#}=(V_{1}^{\#},\ldots,V_{d}^{\#})$.
\end{enumerate}
\end{Alg}
\end{svgraybox}

Repeating the above algorithm $n$ times independently, we get a sample of $n$ independent random vectors
drawn from $\copb$, conditional on the data $\vc{X}_1, \ldots, \vc{X}_n$.  Let this sample be denoted by $\vc{V}_i^{\#} = (V_{i1}^{\#}, \ldots, V_{id}^{\#})$, $i = 1, \ldots, n$.  This procedure can be viewed as a kind of \emph{smoothed bootstrap}
(see \cite{Efron1982siam}, \cite[Section 3.5]{Shao-Tu95}) because the empirical beta copula may be thought of as a smoothed version of the empirical copula.

The joint and marginal empirical distribution functions of the bootstrap sample are
\begin{align*}
  \bbG_n^{\#}( \vc{u} )
  =
  \frac{1}{n} \sum_{i=1}^n \prod_{j=1}^d \I \{ V_{ij}^{\#} \le u_j \} \quad\mbox{and}\quad
  \bbG_{nj}^{\#}(u_j)
  =
  \frac{1}{n} \sum_{i=1}^n \I \{ V_{ij}^{\#} \le u_j \},
\end{align*}
respectively.  The ranks of the bootstrap sample are given by
\begin{equation}
  \label{eq:beta-boot-rank}
  R^{\#}_{ij,n} 
  = n \, \bbG_{nj}^{\#}(V_{ij}^{\#}) 
  = \sum_{k=1}^n \I \{ V_{kj}^{\#} \leq V_{ij}^{\#} \}.
\end{equation}
These yield bootstrapped versions of the Deheuvels empirical copula~\cite{Deheu79}, rank-based empirical copula \eqref{eq:cop:rank}, and empirical beta copula:
\begin{align*}
  \cop^{\#}(\vc{u})
  &:=
  \bbG_n^{\#} \bigl( \bbG_{n1}^{\# -}(u_1), \ldots, \bbG_{nd}^{\# -}(u_d) \bigr),\quad
  \tilbbC_n^{\#}(\vc{u}) 
  := \frac{1}{n} \sum_{i=1}^n \prod_{j=1}^d \I \{ R^{\#}_{ij,n}/n \le u_j \}, \\
  \copbh(\vc{u}) 
  &:= \frac{1}{n} \sum_{i=1}^n \prod_{j=1}^d F_{n,R^{\#}_{ij,n}} (u_j),  
\end{align*}
respectively.  

\begin{Prop}
	\label{prop:beta}
	Assume Condition~\ref{cond:partial-deriv}. Then, as $n \to \infty$, we have conditional weak convergence in probability, as defined in \eqref{eq:weakcondWinP}, with respect to the random vectors $\vecV_1^{\#}, \ldots, \vecV_n^{\#}$ of the bootstrapped empirical copula processes
	\begin{align*}
	\alpha_n^{\#}
	:= \sqrt{n} ( \cop^{\#} - \cop ), \quad
	\tilde{\alpha}_n^{\#}
	:= \sqrt{n} ( \tilbbC_n^{\#} - \tilbbC_n ), \quad
	\alpha_n^{\beta \#}
	:= \sqrt{n} ( \copbh - \copb ),
	\end{align*}
	to the limit process $\bbG^C$ defined in Theorem~\ref{thm:conv-emp-cop}.
\end{Prop}

\subsection{Approximating the sampling distributions of rank statistics by resampling from
  the empirical beta copula}
\label{subsec:beta-resample}
Statistical inference for $C$ often involves rank statistics. 
One way to justify this is to appeal to the invariance of $C$ under coordinatewise, continuous, 
strictly increasing transformations.  Hence we consider a rank statistic $T(\vecR_1,\ldots,\vecR_n)$,
where $\vecR_i:=(R_{i1,n}, \ldots, R_{id,n})$ is a vector of the coordinatewise ranks of $\vc{X}_i$.  
Below, we suggest a way of approximating its distribution by drawing a sample from $\bbC_n^\beta$, and then computing
the ``bootstrap replicates''.  This also avoids problems with the ties encountered when drawing with replacement 
from the original data.  The procedure is as follows.

\begin{svgraybox}\vspace{-\baselineskip}
\begin{Alg}[Smoothed beta bootstrap]
Given $\vecR_1,\ldots,\vecR_n$:
\begin{enumerate}[1.]
\setlength{\itemsep}{5pt}
\item Apply Algorithm~\ref{alg:beta-sampling} $n$ times independently to obtain a bootstrap sample 
$\vc{V}_1^{\#}, \ldots, \vc{V}_n^{\#}$ drawn from $\copb$, compute their ranks $\vc{R}_1^{\#}, 
\ldots, \vc{R}_n^{\#}$ as in \eqref{eq:beta-boot-rank}, and put $T^{\#}:=T(\vc{R}_1^{\#}, \ldots, \vc{R}_n^{\#})$.
\item Repeat Step~1 a moderate-to-large number of times, $B$, to obtain the bootstrap replicates 
$T^{\#}_1, \ldots, T^{\#}_B$.  
\item Use $T^{\#}_1, \ldots, T^{\#}_B$ to approximate the sampling distribution of 
$T(\vecR_1,\ldots,\vecR_n)$.
\end{enumerate}
\end{Alg}
\end{svgraybox}

The validity of this procedure follows from our claim in the preceding subsection.
Because the related empirical copula processes are all asymptotically equivalent,
we need to examine the small-sample performance of the methods.
In Subsection \ref{subsec:rank-corr},  we construct confidence intervals for several copula functionals using popular rank statistics.

\section{Simulation Studies}
\label{sec:sim} 

We assess the performance of the bootstrap methods presented in Sections~\ref{sec:bs-review} and~\ref{sec:beta} in a wide range of applications. In all of the experiments below, the number of Monte Carlo runs and the number of bootstrap replications are both set to $1000$.  We use Clayton, Gumbel--Hougaard, Frank, and Gauss copula families; see, for example, \cite{Nelsen2006}. Most simulations are performed in \textsf{R} using the package \textsf{copula} \cite{copulaR}d; however, the simulation described in Subsection~\ref{subsec:rank-corr} uses MATLAB.  


\subsection{Covariance of the limiting process}

We compare the estimated covariances of the limiting process $\mathbb{G}^C$ based on the standard and smoothed beta bootstrap methods with those of the partial derivatives multiplier method.  In \cite{Buech-Dette2010}, the latter is shown to outperform the straightforward bootstrap and the direct multiplier method.  We follow the setup in \cite{Buech-Dette2010}, evaluating the covariance at four points $\{(i/3,j/3)\}$ for $i,j \in \{1,2\}$ in the unit square. The variables $\xi_1,\ldots,\xi_n$ for the partial derivatives multiplier method are such that $\mathbb{P}[\xi_i = 0] = \mathbb{P} [\xi_i = 2] = 1/2$ for $i \in \{1,\ldots,n\}$. For the bivariate Clayton copula with parameter $\theta = 1$,  Table~\ref{tab:covar} shows the mean squared error of the estimated covariance based on the partial derivatives multiplier method $\alpha_n^{\textnormal{pdm} \circ}$, standard beta bootstrap $\alpha_n^{\beta}$, and smoothed beta bootstrap $\alpha_n^{\beta \#}$, for $n = 100$ and $n = 200$. The results for $\alpha_n^{\textnormal{pdm} \circ}$ are copied from Tables~3 and~4 in \cite{Buech-Dette2010}. Both methods based on the empirical beta copula outperform the multiplier method for all points other than $(1/3,1/3)$ and $(2/3,2/3$).

\begin{table}[ht]
 \caption{Mean squared error $(\times 10^4)$ of the covariance estimates for the bivariate Clayton copula with $\theta = 1$.}  
 \label{tab:covar}
 \begin{center}
 \begin{tabular}{lccccccccc}
\toprule
\multicolumn{2}{l}{} & \multicolumn{4}{c}{$n = 100$ } & \multicolumn{4}{c}{$n = 200$} \\
\cmidrule(r){3-6}
\cmidrule(r){7-10}
\text{} & \text {} & $\left( \tfrac{1}{3},\tfrac{1}{3} \right) \,\,$ & $\left( \tfrac{1}{3},\tfrac{2}{3} \right) \,\,$ & $\left( \tfrac{2}{3},\tfrac{1}{3} \right)\,\,$ & $\left( \tfrac{2}{3},\tfrac{2}{3} \right) \,\,$ & $\left( \tfrac{1}{3},\tfrac{1}{3} \right) \,\,$ & $\left( \tfrac{1}{3},\tfrac{2}{3} \right) \,\,$ & $\left( \tfrac{2}{3},\tfrac{1}{3} \right) \,\,$ & $\left( \tfrac{2}{3},\tfrac{2}{3} \right)$  \\
\midrule
$\alpha_n^{\textnormal{pdm} \circ}$ & $(1/3,1/3)$ &  0.8887 & 0.5210 & 0.5222 & 0.3716 & 0.4595 & 0.2673 & 0.2798 & 0.1961 \\
 & $(1/3,2/3)$ &  & 1.0112 & 0.1799 & 0.2988 &  & 0.5211 & 0.1069 & 0.1577 \\
 & $(2/3,1/3)$ & & & 0.9899 & 0.2818 & & & 0.5092 & 0.1681 \\
 & $(2/3,2/3)$ & & & & 0.6250 & & & & 0.2992 \\
\midrule
$\alpha_n^{\beta}$ & $(1/3,1/3)$ &  0.9992 & 0.3402 & 0.3473 & 0.1956 &  0.6205 & 0.2427 & 0.2383 & 0.1547 \\
 & $(1/3,2/3)$ &  & 0.7887 & 0.1294 & 0.1889 &  & 0.4933 & 0.0857 & 0.1366 \\
 & $(2/3,1/3)$ & & & 0.7644 & 0.1821 & & & 0.4898 & 0.1376 \\
 & $(2/3,2/3)$ & & & & 0.7108 & & & & 0.4183 \\
\midrule
$\alpha_n^{\beta \#}$ & $(1/3,1/3)$ &  1.2248 & 0.2929 & 0.2924 & 0.1456 &  0.6761 & 0.1874 & 0.1888 & 0.1128 \\
 & $(1/3,2/3)$ &  & 0.8461 & 0.0992 & 0.1691  &  & 0.4814 & 0.0703 & 0.1071 \\
 & $(2/3,1/3)$ & & & 0.8856 & 0.1682 & & & 0.4956 & 0.1149 \\
 & $(2/3,2/3)$ & & & & 1.1209 & & & & 0.5913 \\
 \bottomrule
\end{tabular}
  \end{center}
\end{table}

\subsection{Confidence intervals for rank correlation coefficients}
\label{subsec:rank-corr}
Here, we assess the performance of the straightforward bootstrap and the smoothed beta bootstrap (Subsections~\ref{subsec:straight-bs} and \ref{subsec:beta-resample}) for constructing confidence intervals for two popular rank correlation coefficients for bivariate distributions, namely, Kendall's $\tau$ and Spearman's $\rho$, which are known to depend only on the copula $C$ associated with $F$. 

The population Kendall's $\tau$ is defined by
\[
  \tau(C):= 4\int_0^1 \int_0^1 C(u_1,u_2)\,\diff C(u_1,u_2)-1.
\]
In terms of
\begin{gather*}
  Q_{k,i}:=\mathrm{sign}[(X_{k,1}-X_{i,1})(X_{k,2}-X_{i,2})]
  =\mathrm{sign}[(R_{k1,n}-R_{i1,n})(R_{k2,n}-R_{i2,n})], \\
\intertext{and}
  K:=\sum_{i=1}^{n-1}\sum_{k=i+1}^n Q_{k,i},
\end{gather*}
the sample Kendall's $\tau$ is given by $\hat{\tau}:=2K/[n(n-1)]$. Its asymptotic variance can be estimated by
\[
  \hat{\sigma}_\tau^2:=\frac{2}{n(n-1)}\left[\frac{2(n-2)}{n(n-1)^2}
    \sum_{i=1}^n(C_i-\overline{C})^2+1-\hat{\tau}^2\right],
\]
where $C_i:=\sum_{k=1,\ k\not= i}^n Q_{k,i}$, $i\in\{1,\ldots, n\}$ and
$\overline{C}=n^{-1}\sum_{i=1}^n C_i=2K/n$ (see \cite{Hol-Wol-Chi2014}).
Thus, an asymptotic confidence interval for $\tau$ is given by $\hat{\tau}\pm z_{\alpha/2}\hat{\sigma}_\tau$, where
$z_{\alpha/2}$ is the usual standard normal tail quantile. 

This interval can be compared with the confidence intervals obtained using our resampling methods. Table~\ref{tab:conf-tau} shows the coverage probabilities and the average lengths of the estimated confidence intervals based on the asymptotic distribution, straightforward bootstrap, and smoothed beta bootstrap for the independence copula ($\tau = 0$) and the Clayton copula with $\theta = 2$ ($\tau = 0.5$) and $\theta = -2/3$ ($\tau = -0.5)$. The nominal confidence level is 0.95. The smoothed beta bootstrap gives the most conservative coverage probabilities, but has the shortest length of the three. 

\begin{table}[ht]
  \caption{Coverage probabilities and average lengths of the confidence intervals for Kendall's $\tau$ for the Clayton copula family, computed using the normal approximation, straightforward bootstrap, and smoothed beta bootstrap.}
  \label{tab:conf-tau}
  \begin{center}
    \begin{tabular}{llrrrrrrrrrrrr} 
    \toprule
    \multicolumn{2}{l}{} & \multicolumn{4}{c}{$\tau = 0$} & \multicolumn{4}{c}{$\tau = 0.5$} & \multicolumn{4}{c}{$\tau = -0.5$} \\
\cmidrule(r){3-6}
\cmidrule(r){7-10}
\cmidrule(r){11-14}
      & $n$ & 40 & 60 & 80 & 100 & 40 & 60 & 80 & 100 & 40 & 60 & 80 & 100 \\
      \midrule
      coverage & asymp &  0.952 & 0.930 & 0.941 & 0.959 & 0.946 &  0.931 & 0.937 & 0.943 & 0.933 & 0.941 & 0.939 & 0.926 \\
      probability & boot &  0.957 & 0.937 & 0.942 & 0.963 & 0.949 & 0.940 &  0.949 & 0.949 & 0.951 & 0.947 & 0.938 & 0.935 \\
      & beta & 0.964 & 0.949 & 0.949 & 0.966 &0.952 &  0.947 & 0.954 & 0.955 & 0.963 & 0.935 & 0.948 & 0.939 \\
      \midrule
      average & asymp & 0.449 & 0.355 & 0.304 & 0.271 & 0.364 & 0.287 & 0.245 & 0.217 & 0.378 & 0.302 & 0.257 & 0.227 \\
      length & boot & 0.450 & 0.357 &  0.306  & 0.272 & 0.366 & 0.288 & 0.246 & 0.218 & 0.380 & 0.304 & 0.258 & 0.228 \\
      & beta & 0.433 & 0.347 & 0.299  & 0.268 & 0.350 & 0.279 & 0.240 & 0.213 & 0.365 & 0.294 & 0.253 & 0.224 \\ \bottomrule
    \end{tabular}
  \end{center}
\end{table}%

The population Spearman's $\rho$ and the sample Spearman's rho are given by
\begin{align*}
  \rho(C) &:= 12\int_0^1 \int_0^1 \bigl[C(u_1,u_2)-u_1u_2\bigr]\,du_1du_2, \\
  \hat{\rho} & :=\frac{12}{n(n^2-1)}
              \sum_{i=1}^n\left(R_{i1,n}-\frac{n+1}{2}\right)\left(R_{i2,n}-\frac{n+1}{2}\right), 
\end{align*}
respectively.  
The limiting distribution of $\hat{\rho}$ is equal to that of $12\iint\bbG^C(u_1,u_2){\d}u_1{\d}u_2$; thus, in principle, 
it is possible to construct confidence intervals based on the asymptotics.  However, unlike the case of $\hat{\tau}$, this procedure is cumbersome and involves partial derivatives of $C$, which must be estimated.  Therefore, we omit it 
from our study.
We continue to set the nominal confidence level to 0.95 in the experiment. Table~\ref{tab:conf-rho} shows that the coverage probabilities for the smoothed beta bootstrap are more conservative than those for the straightforward bootstrap; however, the average lengths of the estimated confidence intervals are very similar in the two methods.  This could be due to the fact that $\rho(\bbC_n^\beta)=[(n-1)/(n+1)]\hat{\rho}$, as can be computed directly. 

\begin{table}[ht]
  \caption{Coverage probabilities and average lengths of the confidence intervals for Spearman's $\rho$ for the Clayton copula family, based on the straightforward bootstrap and smoothed beta bootstrap.}
  \label{tab:conf-rho}
  \begin{center}
    \begin{tabular}{llrrrrrrrrrrrr} 
    \toprule
    \multicolumn{2}{l}{} & \multicolumn{4}{c}{$\rho = 0$} & \multicolumn{4}{c}{$\rho = 0.5$} & \multicolumn{4}{c}{$\rho = -0.5$} \\
\cmidrule(r){3-6}
\cmidrule(r){7-10}
\cmidrule(r){11-14}
      & $n$ & 40 & 60 & 80 & 100 & 40 & 60 & 80 & 100 & 40 & 60 & 80 & 100 \\
      \midrule
      coverage & boot & 0.956 & 0.943 & 0.953 & 0.951 & 0.959 &0.953 & 0.949 &  0.952 & 0.952 & 0.954 & 0.960 & 0.956\\
      probability & beta & 0.965 & 0.946 & 0.957 & 0.956 & 0.961 & 0.958 & 0.960 & 0.952 & 0.969 & 0.957 & 0.964 & 0.958 \\\midrule
      average & boot & 0.634 & 0.514 & 0.444 & 0.397 & 0.524 & 0.424 & 0.367 & 0.326 & 0.519 & 0.418 & 0.366 & 0.324 \\
      length & beta & 0.625 & 0.510 & 0.442 & 0.395 & 0.522 & 0.424 & 0.368 & 0.325 & 0.519 & 0.418 & 0.367 & 0.324 \\\bottomrule
    \end{tabular}
  \end{center}
\end{table}%

\subsection{Confidence intervals for a copula parameter}

Suppose that the copula of $F$ is parametrized by $\theta\in\Theta\subset\bbR$, such that 
$F(x_1,x_2) = C_\theta(F_1(x_1),F_2(x_2))$. When the $F_j$'s are unknown, the resulting problem of estimating $\theta$
is semiparametric; see \cite{Gen-Gho-Riv95, Tsuka05}.
Assume that $C_\theta$ is absolutely continuous with density $c_\theta$, which is differentiable
with respect to $\theta$.  Replacing the unknown $F_j$'s in the score equation with their (rescaled) empirical 
counterparts, one gets the estimating equation
\begin{equation}
\label{eq:ple}
  \sum_{k=1}^n \frac{\dot{c}_\theta [\bbF_{n1}(X_{k,1}), \bbF_{n2}(X_{k,2})]}
  {c_\theta [\bbF_{n1}(X_{k,1}), \bbF_{n2}(X_{k,2})]} = 0,
\end{equation}
where $\dot{c}_\theta=\partial c_\theta/\partial\theta$.  The solution $\widehat{\theta}$ to \eqref{eq:ple}
is called the \textit{pseudo-likelihood estimator}. 

We compare the confidence intervals for $\theta$ estimated using the pseudo-likelihood estimator $\widehat{\theta}$ based on the asymptotic variance given in \cite{Gen-Gho-Riv95}, straightforward bootstrap, smoothed beta bootstrap, and classic parametric bootstrap.  We set the nominal confidence level equal to 0.95.  Tables~\ref{tab:param-clayton} and~\ref{tab:param-others} show the estimated coverage probabilities and average interval lengths of the confidence intervals for the Clayton, Gauss, Frank, and Gumbel--Hougaard copula families, respectively.  For the Clayton copula, the smoothed beta bootstrap gives the shortest intervals, both for $\theta=1$ and $\theta=2$, but for $\theta=2$, the coverage probabilities are too liberal, which is somewhat puzzling.  For the Frank and Gumbel--Hougaard copulas, the smoothed beta bootstrap gives the most conservative coverage probabilities, but has the shortest length of the four.  For the Gauss copula, the asymptotic approximation gives significantly smaller coverage probabilities than the nominal value of 0.95. 

\begin{table}[ht]
   \caption{Coverage probabilities and average lengths of the confidence intervals for the parameter of the Clayton copula, with $\theta =1$ ($\tau = 1/3$) and $\theta = 2$ ($\tau = 1/2$). Intervals are computed using the asymptotic normal approximation, straightforward bootstrap, smoothed beta bootstrap, and parametric bootstrap.}
  \label{tab:param-clayton}
  \begin{center}
    \begin{tabular}{llrrrrrrrr} 
    \toprule
    \multicolumn{2}{l}{} & \multicolumn{4}{c}{$\theta = 1$} & \multicolumn{4}{c}{$\theta = 2$} \\
\cmidrule(r){3-6}
\cmidrule(r){7-10}
      & $n$ & 40 & 60 & 80 & 100 & 40 & 60 & 80 & 100 \\
      \midrule
      coverage & asymp & 0.954 & 0.969 & 0.960 & 0.965 & 0.951 & 0.940 & 0.940 & 0.946 \\
      probability & boot & 0.953 & 0.943 & 0.944 & 0.943 & 0.968 & 0.952 & 0.953 &  0.951 \\
      & beta & 0.953 & 0.964 & 0.957 & 0.952 & 0.933 & 0.904 & 0.908 & 0.906 \\  
      & param & 0.924 & 0.923 & 0.933  & 0.948 & 0.957 & 0.951 & 0.955 & 0.953
       \\\midrule
      average & asymp & 2.011 & 1.632 & 1.354 & 1.237 & 2.764 & 2.142 & 1.821 & 1.615 \\
      length & boot & 1.894 & 1.449 & 1.198 & 1.046 & 2.991 & 2.205 & 1.841 & 1.626 \\ 
      & beta & 1.517 & 1.225 & 1.050 & 0.935 & 1.957 & 1.612 & 1.420 & 1.296 \\
       & param & 1.914 & 1.448 & 1.222 & 1.070 & 2.821 & 2.150 & 1.829 & 1.617 \\
        \bottomrule
    \end{tabular}
  \end{center}
\end{table}%

\begin{table}[ht]
   \caption{Coverage probabilities and average lengths of confidence intervals for the parameter of the Gaussian copula with $\theta = 1/\sqrt{2}$, the Frank copula with $\theta = 5.75$ and the Gumbel--Hougaard copula with $\theta =  2$. All copulas have $\tau \approx 1/2$. Intervals computed via the asymptotic normal approximation, the straightforward bootstrap, the smoothed beta bootstrap, and the parametric bootstrap.}
  \label{tab:param-others}
  \begin{center}
    \begin{tabular}{llrrrrrrrrrrrr} 
    \toprule
    \multicolumn{2}{l}{} & \multicolumn{4}{c}{Gauss} & \multicolumn{4}{c}{Frank} & \multicolumn{4}{c}{Gumbel--Hougaard} \\
\cmidrule(r){3-6}
\cmidrule(r){7-10}
\cmidrule(r){11-14}
      & $n$ & 40 & 60 & 80 & 100 & 40 & 60 & 80 & 100 & 40 & 60 & 80 & 100 \\
      \midrule
      coverage & asymp & 0.881 & 0.895 & 0.910 & 0.928 & 0.941 & 0.950 & 0.948 & 0.965 & 0.954 & 0.940 & 0.940 & 0.955 \\
      probability & boot & 0.942 & 0.944 & 0.947 & 0.950 & 0.957 & 0.956 & 0.946 & 0.963 & 0.965 & 0.951 & 0.953 & 0.965 \\
      & beta & 0.968 & 0.962 & 0.970 & 0.953 & 0.965 & 0.961 & 0.952 & 0.965 & 0.970 & 0.951 & 0.952 & 0.954 \\
      & param &  0.903 & 0.921 & 0.923 & 0.930 & 0.938 & 0.956 & 0.941 & 0.962 & 0.924 & 0.926 & 0.932 & 0.945 \\
       \midrule
      average & asymp & 0.303 & 0.274 & 0.213 & 0.193 & 5.699 & 4.487 & 3.821 & 3.391 & 1.425 & 1.082 & 0.929 & 0.816 \\
      length & boot & 0.319 & 0.257 & 0.219 & 0.197 & 6.139 & 4.677 & 3.949 & 3.464 & 1.572 & 1.162 & 0.968 & 0.855 \\
      & beta & 0.341 & 0.269 & 0.228 & 0.203 & 5.367 & 4.335 & 3.735 & 3.329 & 1.170 & 0.947 & 0.826 & 0.747 \\
       & param & 0.292 & 0.242 & 0.210 & 0.191 & 5.729 & 4.494 & 3.848 & 3.389 & 1.546 & 1.170 & 0.983 & 0.869 \\
        \bottomrule
    \end{tabular}
  \end{center}
\end{table}%

\subsection{Testing the symmetry of a copula}

For a bivariate copula $C$, consider the problem of testing the symmetry hypothesis $H_0: C(u_1, u_2)  = C(u_2,u_1)$ for all $(u_1, u_2) \in [0,1]^2$. We focus on the following two test statistics proposed in \cite{Gen-Nes-Que2012}:
\begin{align*}
  S_n & =\int_{[0,1]^2}\left[\bbC_n(u_1,u_2)-\bbC_n(u_2,u_1)\right]^2\,{\d}\bbC_n(u_1,u_2), \\
  R_n & =  \int_{[0,1]^2} \left[ \cop (u_1,u_2) - \cop (u_2,u_1) \right]^2 \,{\d} u_1 \,{\d} u_2,
\end{align*}
and include versions based on the empirical beta copula; that is,
\begin{align*}
  S_n^\beta & =\int_{[0,1]^2}\left[\bbC_n^\beta(u_1,u_2)-\bbC_n^\beta(u_2,u_1)\right]^2\,{\d}\bbC_n^\beta(u_1,u_2), \\
	R_n^\beta &  =
	\int_{[0,1]^2} [ \bbC_n^\beta(u_1, u_2) - \bbC_n^\beta(u_2, u_1) ]^2 \, 
	\diff u_1 \, \diff u_2.
\end{align*}
Similarly, as in Proposition~1 in \cite{Gen-Nes-Que2012}, the statistic $R_n^\beta$ can be computed as
\begin{multline*}
	R_n^\beta
	= \frac{2}{n^2} \sum_{i=1}^n \sum_{j=1}^n
	\{ B_n(R_{i1,n}, R_{j1,n}) B_n(R_{i2,n}, R_{j2,n})
        - \\
	B_n(R_{i1,n}, R_{j2,n}) B_n(R_{i2,n}, R_{j1,n}) \},
\end{multline*}
with $B_n(r, s) = \int_0^1 F_{n,r}(u) F_{n,s}(u) \, \diff u$ for $r,s \in \{1, \ldots, n\}$ and $F_{n,r}(u)$ as in \eqref{eq:betacdf}. For fixed $n$, the matrix $B_n$ can be precomputed and stored, which reduces the computation time for the resampling methods. Similarly, $S_n^\beta$ can be written as
\begin{align*}
	S_n^\beta
	= n^{-3}  \sum_{i=1}^n \sum_{j=1}^n \sum_{k=1}^n
	\big\{ & C_n(R_{i1,n}, R_{j1,n} R_{k1,n}) C_n(R_{i2,n}, R_{j2,n} R_{k2,n})
        \\
        &  \, - C_n(R_{i1,n}, R_{j2,n} R_{k1,n}) C_n(R_{i2,n}, R_{j1,n} R_{k2,n}) \\
             & \, - C_n(R_{i2,n}, R_{j1,n} R_{k1,n}) C_n(R_{i1,n}, R_{j2,n} R_{k2,n}) \\
                     & \, + C_n(R_{i2,n}, R_{j2,n} R_{k1,n}) C_n(R_{i1,n}, R_{j1,n} R_{k2,n})    	 \big\},
\end{align*}
with $C_n(r, s,t) = \int_0^1 F_{n,r}(u) F_{n,s}(u) \, \diff F_{n,t} (u)$ for $r,s,t \in \{1, \ldots, n\}$.

In order to compute the $p$-values, we need to generate bootstrap samples from a distribution that fulfills the restriction specified by $H_0$. A natural candidate is a ``symmetrized'' version of the empirical beta copula
\[
  \bbC_n^{\beta,\mathrm{sym}}(u_1,u_2)
  :=
  \half\bbC_n^\beta(u_1,u_2)+\half\bbC_n^\beta(u_2,u_1).
\]
When resampling, this simply amounts to interchanging the two coordinates at random in step 3 of Algorithm~\ref{alg:beta-sampling}.  We employ the following three resampling schemes to compare the actual sizes of the tests.
\begin{itemize}
\item The symmetrized smoothed beta bootstrap: we resample from $\bbC_n^{\beta,\mathrm{sym}}$ to obtain bootstrap replicates of $R_n$, $R_n^\beta$, $S_n$, and $S_n^{\beta}$;
\item The symmetrized version of the straightforward bootstrap for $R_n$ and $S_n$;
\item \textsl{exchTest} in the \textsf{R} package \textsf{copula} \cite{copulaR}, which implements the multiplier bootstrap for $R_n$ and $S_n$, as described in \cite{Gen-Nes-Que2012} and in Section~5 of \cite{Kojadinovic-Yan2012}. For $R_n$, the grid length in \textsl{exchTest} is set to $m = 50$.
\end{itemize}

We use the nominal size $\alpha=0.05$ throughout the experiment. Tables~\ref{tab:sym-clayton} and \ref{tab:sym-gauss} show the actual sizes of the symmetry tests for the Clayton and Gauss copulas.  On the whole, the smoothed beta bootstrap works better than \textsl{exchTest}, and works equally well as $R_n$ and $S_n$, except when the dependence is strong ($\tau=0.75$) and the sample size is small ($n = 50$).  However, no method produces a satisfying result in the latter case.  The smoothed beta bootstraps with $R_n^\beta$ and $S_n^\beta$ produce actual sizes similar to those based on $R_n$ and $S_n$.  The statistic $S_n$ performs slightly better than $R_n$ on average, especially for strong positive dependence. The straightforward bootstrap performs poorly in all cases, as expected \cite{Remillard-Scaillet2009}.

To compare the power of the tests, the Clayton and Gauss copulas are made asymmetric using Khoudraji's device~\cite{Khoudraji1995}; that is, the asymmetric version of a copula $C$ is defined as
\begin{equation*}
	K_{\delta} (u_1, u_2) = u_1^{\delta} C(u_1^{1-\delta}, u_2), 
	\qquad (u_1, u_2) \in [0,1]^2.
\end{equation*}
Table \ref{tab:sym-power} shows the empirical power of $R_n$ and $R_n^\beta$ for $\delta \in \{0.25,0.5,0.75\}$ for the three resampling methods.  As shown, the smoothed beta bootstraps with $R_n$ and $R_n^\beta$ have higher power than \textsl{exchTest} for almost all sample sizes and parameter values considered; furthermore, the smoothed beta bootstrap with $R_n^\beta$ has slightly higher power in almost all cases.

\begin{table}[ht]
  \caption{Actual sizes of the symmetry tests based on $R_n$ and $S_n$ for the Clayton copula ($\theta \in \{-1/3, 2/3, 2, 6 \}$), with $p$-values computed by the multiplier bootstrap (exchTest), straightforward bootstrap (boot), smoothed beta bootstrap (beta), and of the test based on $R_n^\beta$ and $S_n^\beta$, with $p$-values computed by the smoothed beta bootstrap (beta2). The nominal size is $\alpha=0.05$}
  \label{tab:sym-clayton}
  \begin{center}
    \begin{tabular}{ccrrrrcccrrrr} \toprule
      \multicolumn{1}{c}{} &   \multicolumn{5}{c}{$S_n$} &      \multicolumn{1}{c}{} &  \multicolumn{5}{c}{$R_n$}  \\
\cmidrule(r){2-6}
\cmidrule(r){8-12}
      & $n$ & 50 & 100  & 200 & 400 & &  $n$ & 50 & 100  & 200 & 400  \\\midrule
      $\tau=-1/5$ & \textsf{exchTest}  & 0.055 & 0.033 & 0.039 & 0.040 &       $\tau=-1/5$ & \textsf{exchTest} & 0.044 & 0.035 & 0.039 & 0.051 \\
 & boot & 0.021 & 0.024 & 0.031 & 0.034 &  & boot & 0.009 & 0.019 & 0.027 & 0.044 \\
      & beta & 0.057 & 0.038 & 0.039 & 0.041 & & beta & 0.046 & 0.035 & 0.042 & 0.059 \\
      & beta2 & 0.042 & 0.036 & 0.043 & 0.057 & & beta2 & 0.050 & 0.037 & 0.041 & 0.059 \\
      \midrule
$\tau=0.25$& \textsf{exchTest} & 0.039 & 0.029 & 0.036 & 0.036 & $\tau=0.25$ & \textsf{exchTest} & 0.030 & 0.022 & 0.040 & 0.031 \\
      & boot  & 0.009 & 0.015 & 0.026 & 0.034 &    & boot & 0.001 & 0.011 & 0.020 & 0.030 \\
      & beta & 0.039 & 0.039 & 0.044 & 0.046 &   & beta & 0.042 & 0.032 & 0.041 & 0.043 \\ 
      & beta2 & 0.043 & 0.034 & 0.041 & 0.044 & & beta2 & 0.033 & 0.033 & 0.044 & 0.045 \\
      \midrule
     $\tau=0.5$ & \textsf{exchTest} & 0.033 & 0.020 & 0.026 & 0.019 &  $\tau=0.5$ & \textsf{exchTest} & 0.015 & 0.014 & 0.030 & 0.031 \\
      & boot & 0.001 & 0.008 & 0.015 & 0.017 &  & boot &  0.001 & 0.005 & 0.017 & 0.025 \\	
      & beta & 0.030 & 0.029 & 0.039 & 0.028 &  & beta & 0.020 & 0.022 & 0.040 & 0.047 \\ 
      & beta2 & 0.029 & 0.024 & 0.046 & 0.041 & & beta2 & 0.019 & 0.023 & 0.045 & 0.046 \\
      \midrule
   $\tau=0.75$& \textsf{exchTest} & 0.025 & 0.026 & 0.018 & 0.014 &   $\tau=0.75$& \textsf{exchTest} & 0.000 & 0.007 & 0.007  & 0.012 \\
     & boot & 0.000 & 0.002 & 0.001 & 0.008 &  & boot & 0.000 & 0.000 & 0.001 & 0.004 \\
      & beta  & 0.006 & 0.017 & 0.026 & 0.029 &  & beta & 0.000 & 0.006 & 0.017 & 0.029 \\
      & beta2 & 0.001 & 0.014 & 0.029 & 0.034 & & beta2 & 0.000 & 0.007 & 0.021 & 0.036 \\
      \bottomrule
    \end{tabular}
  \end{center}
\end{table}%

\begin{table}[ht]
  \caption{Actual sizes of the symmetry tests based on $R_n$ and $S_n$ for the Gauss copula, with $p$-values computed by the multiplier bootstrap (exchTest), straightforward bootstrap (boot), smoothed beta bootstrap (beta), and of the test based on $R_n^\beta$ and $S_n^\beta$, with $p$-values computed by the smoothed beta bootstrap (beta2).  The nominal size is $\alpha=0.05$}
  \label{tab:sym-gauss}
  \begin{center}
    \begin{tabular}{ccrrrrcccrrrr} \toprule
      \multicolumn{1}{c}{} &   \multicolumn{5}{c}{$S_n$} &      \multicolumn{1}{c}{} &  \multicolumn{5}{c}{$R_n$}  \\
\cmidrule(r){2-6}
\cmidrule(r){8-12}
      & $n$ & 50 & 100  & 200 & 400 & &  $n$ & 50 & 100  & 200 & 400  \\\midrule
$\tau=-0.5$  & \textsf{exchTest} & 0.047 & 0.032 & 0.038 & 0.039 & $\tau=-0.5$  & \textsf{exchTest} & 0.026 & 0.037 & 0.037 & 0.041 \\
       & boot & 0.022 & 0.023 & 0.032 & 0.040 & & boot & 0.007 & 0.014 & 0.027 & 0.033 \\ 
       & beta & 0.044 & 0.030 & 0.035 & 0.043 & & beta & 0.020 & 0.030 & 0.036 & 0.040 \\
       & beta2 & 0.028 & 0.035 & 0.045 & 0.043 & & beta2 & 0.022 & 0.034 & 0.041 & 0.042 \\
      \midrule
   $\tau=0.25$ & \textsf{exchTest} & 0.028 & 0.035 & 0.031 & 0.040 &  $\tau=0.25$ & \textsf{exchTest} & 0.029 & 0.025 & 0.030 & 0.041 \\
     & boot  & 0.007 & 0.015 & 0.023 & 0.038 &   & boot & 0.008 & 0.015 & 0.025 & 0.037 \\
      & beta & 0.033 & 0.038 & 0.039 & 0.045 &   & beta & 0.040 & 0.030 & 0.033 & 0.048 \\
      & beta2 & 0.048 & 0.033 & 0.037 & 0.047 & & beta2 & 0.037 & 0.031 & 0.037 & 0.048  \\
      \midrule
     $\tau=0.5$ & \textsf{exchTest} & 0.034 & 0.018 & 0.025 & 0.029 & $\tau=0.5$ & \textsf{exchTest} & 0.011 & 0.014 & 0.019 & 0.034 \\
       & boot & 0.003 & 0.005 & 0.015 & 0.026 &  & boot & 0.001 & 0.005 & 0.006 & 0.028 \\
      & beta & 0.032 & 0.033 & 0.041 & 0.044 &  & beta & 0.016 & 0.024 & 0.031 & 0.042 \\
      & beta2 & 0.029 & 0.026 & 0.034 & 0.048 & & beta2 & 0.018 & 0.023 & 0.033 & 0.047 \\
      \midrule
        $\tau=0.75$ & \textsf{exchTest} & 0.018 & 0.017 & 0.011 & 0.008 &       $\tau=0.75$& \textsf{exchTest} & 0.001 & 0.001 & 0.005 &  0.009\\
      & boot & 0.000 & 0.000 & 0.002 & 0.006 & & boot & 0.000 & 0.000 & 0.000 &  0.003 \\
      & beta & 0.006 & 0.015 & 0.021 & 0.029 &  & beta & 0.001 & 0.001 & 0.010 & 0.028 \\
      & beta2 & 0.002 & 0.005 & 0.016  & 0.035 & & beta2 & 0.001 & 0.001 & 0.011 & 0.027 \\
      \bottomrule
    \end{tabular}
  \end{center}
\end{table}%

\begin{table}[h!]
	\caption{Actual power of the symmetry tests based on $R_n$, with $p$-values computed by the multiplier bootstrap (exchTest), straightforward bootstrap (boot), smoothed beta bootstrap (beta), and of the test based on $R_n^\beta$, with $p$-values computed by the smoothed beta bootstrap (beta2), for the Clayton and Gauss copulas, made asymmetric by Khoudraji's device. The nominal size is $\alpha=0.05$}
	\label{tab:sym-power}
	\begin{center}
		\begin{tabular}{ccrrrrcccrrrr} \toprule
			\multicolumn{1}{c}{} &   \multicolumn{5}{c}{Clayton} &      \multicolumn{1}{c}{} &  \multicolumn{5}{c}{Gauss}  \\
			\cmidrule(r){2-6}
			\cmidrule(r){8-12}
			$\delta = 0.25$      & $n$ & 50 & 100  & 200 & 400 & &  $n$ & 50 & 100  & 200 & 400  \\\midrule
			$\tau=0.25$& \textsf{exchTest} & 0.025 & 0.031 & 0.042  & 0.049 & $\tau=0.25$ & \textsf{exchTest} & 0.033 & 0.028 &  0.034 & 0.050  \\
			& boot & 0.006 & 0.017 & 0.029 & 0.042 & & boot & 0.006 & 0.010 & 0.032 & 0.047 \\
                        & beta & 0.034 & 0.039 & 0.056 & 0.050 & & beta & 0.044 & 0.042 & 0.048 & 0.057 \\
			& beta2 & 0.034 & 0.041 &  0.055 & 0.054 &  & beta2 & 0.046 & 0.033 & 0.048  & 0.059  \\
			\midrule
			$\tau=0.5$ & \textsf{exchTest} & 0.047& 0.090 & 0.197  & 0.449  &  $\tau=0.5$ & \textsf{exchTest} & 0.052 & 0.078 & 0.188 & 0.401  \\
			& boot & 0.004 & 0.041 & 0.145  & 0.407 &  & boot & 0.007 & 0.028 & 0.136  & 0.366 \\
                        & beta & 0.060 & 0.100 & 0.216 & 0.469 & & beta & 0.061 & 0.088 & 0.198 & 0.433 \\
			& beta2 & 0.062 & 0.103 & 0.227  & 0.486 &  & beta2 & 0.065 & 0.098  & 0.212  & 0.441  \\
			\midrule
			$\tau=0.75$& \textsf{exchTest} & 0.199 & 0.630  & 0.985  & 1.000 & $\tau=0.75$& \textsf{exchTest} & 0.205 & 0.637 & 0.973  & 1.000 \\
			& boot & 0.038 & 0.380 & 0.949 & 1.000 & & boot & 0.051 & 0.338 & 0.921 & 1.000 \\
                        & beta & 0.227 & 0.637 & 0.981 & 1.000 & & beta & 0.208 & 0.614 & 0.974 & 1.000 \\
			& beta2  & 0.242  &  0.667&  0.986 & 1.000 & & beta2 & 0.225 & 0.639  & 0.986  & 1.000 \\
			\midrule

			$\delta = 0.5$      & $n$ & 50 & 100  & 200 & 400 & &  $n$ & 50 & 100  & 200 & 400  \\\midrule
			$\tau=0.25$& \textsf{exchTest} & 0.028 &  0.029 &  0.050 &  0.055 & $\tau=0.25$ & \textsf{exchTest} & 0.044 & 0.053 & 0.053  & 0.083  \\
			& boot & 0.008 & 0.011 & 0.031 & 0.054 & & boot & 0.009 & 0.019 & 0.034 & 0.068 \\
                        & beta & 0.035 & 0.033 & 0.056 & 0.064 & & beta & 0.051 & 0.059 & 0.055 & 0.090 \\
			& beta2 & 0.038 & 0.037  & 0.059  & 0.064 &  & beta2 & 0.052 & 0.062  & 0.056  & 0.091 \\
			\midrule
			$\tau=0.5$ & \textsf{exchTest} & 0.069 & 0.127 & 0.269 & 0.576 &  $\tau=0.5$ & \textsf{exchTest} & 0.100 & 0.203 & 0.388   & 0.730  \\
			& boot & 0.015 & 0.068 & 0.219 & 0.539 & & boot & 0.027 & 0.119 & 0.326 & 0.695 \\
                        & beta & 0.077 & 0.140 & 0.299 & 0.593 & & beta & 0.105 & 0.213 & 0.410 & 0.741 \\
			& beta2 & 0.075 & 0.144 & 0.306 & 0.602 &  & beta2 & 0.116& 0.225  & 0.416  & 0.750 \\
			\midrule
			$\tau=0.75$& \textsf{exchTest} & 0.385 & 0.814 &  0.997 &  1.000 & $\tau=0.75$& \textsf{exchTest} & 0.478 & 0.914 & 1.000 & 1.000 \\
			& boot & 0.125 & 0.644 & 0.993 & 1.000 & & boot & 0.198 & 0.792 & 1.000 & 1.000 \\
                        & beta & 0.393 & 0.824 & 0.997 & 1.000 & & beta & 0.475 & 0.916 & 1.000 & 1.000 \\
			& beta2  & 0.425 & 0.833 & 0.998 & 1.000 & & beta2 & 0.507 & 0.923 & 1.000 & 1.000 \\
			\midrule

			$\delta = 0.75$      & $n$ & 50 & 100  & 200 & 400 & &  $n$ & 50 & 100  & 200 & 400  \\\midrule
			$\tau=0.25$& \textsf{exchTest} & 0.032& 0.039 & 0.046 & 0.054& $\tau=0.25$ & \textsf{exchTest} & 0.043 & 0.046  & 0.056  & 0.076  \\
			& boot & 0.008 & 0.020 & 0.030 & 0.047 & & boot & 0.016 & 0.027 & 0.036 & 0.060 \\
                        & beta & 0.036 & 0.043 & 0.051 & 0.060 & & beta & 0.047 & 0.053 & 0.061 & 0.081 \\ 
			& beta2 & 0.039 & 0.045 & 0.054   & 0.058 &  & beta2 & 0.053 & 0.053 & 0.060  & 0.081  \\
			\midrule
			$\tau=0.5$ & \textsf{exchTest} & 0.042 & 0.089 & 0.129 & 0.266 &  $\tau=0.5$ & \textsf{exchTest} & 0.088 & 0.169  & 0.317  & 0.655  \\
			& boot & 0.012 & 0.045 & 0.100 & 0.251 & & boot & 0.030 & 0.113 & 0.271 & 0.604 \\
                        & beta & 0.053 & 0.094 & 0.132 & 0.279 & & beta & 0.099 & 0.184 & 0.342 & 0.660 \\
			& beta2 & 0.053 & 0.102 &  0.140 & 0.287 &  & beta2 & 0.102 & 0.194 & 0.354   & 0.670 \\
			\midrule
			$\tau=0.75$& \textsf{exchTest} & 0.144 & 0.372  & 0.693  & 0.962 & $\tau=0.75$& \textsf{exchTest} & 0.369& 0.636 & 0.947 & 1.000 \\
			& boot & 0.051 & 0.268 & 0.645 & 0.958 & & boot & 0.133 & 0.510 & 0.918 & 1.000 \\
                        & beta & 0.156 & 0.382 & 0.720 & 0.965 & & beta & 0.303 & 0.637 & 0.950 & 1.000 \\
			& beta2  & 0.169 &  0.402& 0.727 & 0.966 & & beta2 & 0.334 & 0.664  & 0.954  &  1.000\\
			\bottomrule
		\end{tabular}
	\end{center}
\end{table}%

\section{Concluding Remarks}
\label{sec:concl-rmks}
We have studied the performance of resampling procedures based on the empirical beta copula, and proved that all related empirical copula processes exhibit asymptotically equivalent behavior.  A comparative analysis based on Monte Carlo experiments shows that, on the whole, the smoothed beta bootstrap works fairly well, providing a useful alternative to existing resampling schemes.  However, we also find that its effectiveness varies somewhat between copulas.  

Higher-order asymptotics for the various nonparametric copula estimators might improve our understanding of the various resampling procedures \cite{Hall92, Shao-Tu95}, although calculating such expansions seems a formidable task.  


\begin{acknowledgement}
The research of H. Tsukahara was supported by JSPS KAKENHI Grant Number 18H00836. The authors wish to thank an anonymous reviewer for his/her careful reading of the manuscript and helpful comments.
\end{acknowledgement}

\section*{Appendix: Mathematical Proofs}
\subparagraph{Proof of Proposition~\ref{prop:bs-emp-cop-diff}}
Let $\mathcal{N}_n = \{ i = 1, \ldots, n : W_{ni} \ge 1 \}$ be the set of indices that are sampled at least once. Then, $\bbF_{nj}^*$ is a discrete distribution function with atoms $\{ X_{ij}\colon i \in \mathcal{N}_n \}$ and probabilities $n^{-1} W_{ni}$.

Since $n^{-1} R_{ij,n}^* = \bbF_{nj}^*(X_{ij})$, we have
\begin{align*}
  \left\lvert
    \bbC_n^{*} (\vc{u}) - \tilbbC_n^{*} (\vc{u})
  \right\rvert
  &\le
  \frac{1}{n} \sum_{i\in\mathcal{N}_n} W_{ni}
  \left\lvert
    \prod_{j=1}^d \I \{ X_{ij} \le \bbF_{nj}^{* -}(u_j) \}
    -
    \prod_{j=1}^d \I \{ \bbF_{nj}^*(X_{ij}) \le u_j \}
  \right\rvert \\
  &\le
  \frac{1}{n} \sum_{i\in\mathcal{N}_n} W_{ni}
  \sum_{j=1}^d
    \left\lvert
      \I \{ X_{ij} \le \bbF_{nj}^{* -}(u_j) \}
      -
      \I \{ \bbF_{nj}^*(X_{ij}) \le u_j \}
    \right\rvert \\
  &=
  \frac{1}{n} \sum_{i\in\mathcal{N}_n} W_{ni}
  \sum_{j=1}^d
    \left\lvert
      \I \{ X_{ij} = \bbF_{nj}^{* -}(u_j) \}
      -
      \I \{ \bbF_{nj}^*(X_{ij}) = u_j \}
    \right\rvert.
\end{align*}
In the last equality, we use the fact that $x < G^-(u)$ if and only if $G(x) < u$ for any (right-continuous) distribution function $G$, any real $x$, and any $u \in [0, 1]$. 
For each $j \in \{1, \ldots, d\}$, $\bbF_{nj}^*(X_{ij}) = u_j$ implies $X_{ij} = \bbF_{nj}^{* -}(u_j)$ since $\bbF_{nj}^*$ jumps at $X_{ij}$, $i\in\mathcal{N}_n$, and there is at most a single $i \in \mathcal{N}_n$ such that $X_{ij} = \bbF_{nj}^{* -}(u_j)$.  Thus, we have
\[
  \left\lvert
    \bbC_n^{*} (\vc{u}) - \tilbbC_n^{*} (\vc{u})
  \right\rvert
  \leq
  \frac{d}{n} \max_{i = 1, \ldots, n} W_{ni}. 
\]

By coupling the multinomial random vector $(W_{n1}, \ldots, W_{nn})$ to a vector of independent Poisson(1) random variables $(W_{n1}', \ldots, W_{nn}')$, as in \cite[pp.~346--348]{Vaart-Wellner}, it can be shown that $\max_{i=1,\ldots,n} W_{ni} = O_p(\log n)$ as $n \to \infty$. Equation~\eqref{eq:copEquivCop} follows.


\subparagraph{Proof of Proposition~\ref{prop:asymp-equiv}}
Fix $\eps > 0$.  We know from \eqref{eq:alphatildeconv} 
that $\tilde{\alpha}_n$ converges weakly in $\ell^\infty([0, 1]^d)$ to a Gaussian process with continuous trajectories. Write $\vc{S} = (S_1, \ldots, S_d)$ and for a point $\vc{x} \in \bbR^d$, put $\lvert \vc{x} \rvert_\infty = \max( \lvert x_1 \rvert, \ldots, \lvert x_d \rvert)$. Furthermore, put $\lVert f \rVert_\infty = \sup \{ \lvert f(\vc{u}) \rvert : \vc{u} \in [0, 1]^d \}$ for $f : [0, 1]^d \to \bbR$. Then,
\begin{align*}
  \lvert \alpha_n^\beta(\vc{u}) - \tilde{\alpha}_n(\vc{u}) \rvert
  &\le
  \E_S \left[
    \lvert
      \tilde{\alpha}_n(S_1/n, \ldots, S_d/n) - \tilde{\alpha}_n(\vc{u})
    \rvert
  \right] \\
  &\le
  2 \lVert \tilde{\alpha}_n \rVert_\infty \, \Pr_S [ \lvert \vc{S}/n - \vc{u} \rvert_\infty > \eps ]
  +
  \sup_{\substack{\vc{v}, \vc{w} \in [0, 1]^d \\ \lvert \vc{v} - \vc{w} \rvert_\infty \le \eps}}
  \lvert \tilde{\alpha}_n(\vc{v}) - \tilde{\alpha}_n(\vc{w}) \rvert,
\end{align*}
where $\E_S$ and $\Pr_S$ denote the expectation and probability, respectively, with respect to $\vc{S}$, conditional on 
the sample and the multinomial random vector. 
Let $Y_n(\eps)$ denote the supremum on the right-hand side.  By Tchebysheff's inequality, the probability in the first term on the right-hand side is bounded by a constant multiple of $n^{-1/2} \eps^{-1}$ and thus tends to zero uniformly in $\vc{u} \in [0, 1]^d$. Since $\lVert \tilde{\alpha}_n \rVert_\infty = O_p(1)$, we get
$\lVert \alpha_n^\beta - \tilde{\alpha}_n \rVert_\infty = o_p(1) + Y_n(\eps)$ as $n \to \infty$.
By the weak convergence of $\tilde{\alpha}_n$ in $\ell^\infty([0, 1])$ to a process with continuous trajectories, we can find, for any $\eta > 0$, a sufficiently small $\eps > 0$ such that $\limsup_{n \to \infty} \Pr[Y_n(\eps) > \eta] \le \eta$.
Equation~\eqref{eq:alphabeta:conv} follows. \qed

\subparagraph{Proof of Proposition~\ref{prop:beta}}
\emph{Step 1.} Recall the $C$-pinned Brownian sheet $\mathbb{U}^C$ defined prior to Theorem~\ref{thm:conv-emp-cop}. We show that
\begin{equation}
\label{eq:toshow:gamma}
\gamma_n^{\#} := \sqrt{n} (\bbG_n^{\#} - \copb)
\weakcondVinP \bbU^C, \qquad n \tendsto \infty.
\end{equation}
From~\eqref{eq:weakcondWinP:2}, two claims need to be shown: convergence in the bounded Lipschitz metric (Step~1.1), and asymptotic measurability (Step~1.2).
\smallskip

\emph{Step 1.1} Let $\mathcal{P}$ denote the set of all Borel probability measures on $[0, 1]^d$. For $P \in \mathcal{P}$, let $\bbU^P$ denote a tight, $P$-Brownian bridge on $[0, 1]^d$. Specifically, $\bbU^P$ is a centered Gaussian process with covariance function $\E[ \bbU^P(\vc{u}) \bbU^P(\vc{v})] = F(\vc{u} \wedge \vc{v}) - F(\vc{u}) F(\vc{v})$, where $F$ is the cumulative distribution function associated with $P$, and whose trajectories are uniformly continuous, almost surely with respect to the standard deviation semimetric \cite[Example~1.5.10]{Vaart-Wellner}:
\begin{equation}
\label{eq:sdmetric}
d_P(\vc{u},\vc{v}) 
= \bigl(\E[\{\bbU^P(\vecu) - \bbU^P(\vecv)\}^2]\bigr)^{1/2}
= \bigl( F(\vecu) - 2F(\vecu \wedge \vecv) + F(\vecv) \bigr)^{1/2}.
\end{equation}

Furthermore, for $P \in \mathcal{P}$, let $\bbU_{n,P}$ denote the empirical process based on independent random sampling from $P$. We view $\bbU_{n,P}$ as a random element of $\ell^\infty([0, 1]^d)$ from the empirical and true cumulative distribution functions. Let $P_n^\beta \in \mathcal{P}$ be the (random) probability measure associated with the empirical beta copula $\bbC_n^\beta$. Recall $\operatorname{BL}_1$ in \eqref{eq:BL1}. 

We need to show that, as $n \to \infty$,
\[
\sup_{h \in \operatorname{BL}_1}
\left\lvert
\E_{P_n^\beta}^*[h(\bbU_{n,P_{n}^{\beta}})] - \E[h(\bbU^C)]
\right\rvert
\longrightarrow 0
\qquad \text{in outer probability.}
\]
By the triangle inequality, it is sufficient to show the pair of convergences
\begin{align}
\label{eq:1.1.1}
\sup_{h \in \operatorname{BL}_1}
\left\lvert
\E_{P_n^\beta}^*[h(\bbU_{n,P_n^\beta})] - \E[h(\bbU^{P_{n}^{\beta}})]
\right\rvert
\longrightarrow 0, \\
\label{eq:1.1.2}
\sup_{h \in \operatorname{BL}_1}
\left\lvert
\E[h(\bbU^{P_{n}^{\beta}})] - \E[h(\bbU^C)]
\right\rvert
\longrightarrow 0,
\end{align}
as $n \to \infty$, in outer probability. We do so in Steps~1.1.1 and~1.1.2, respectively.
\smallskip

\emph{Step~1.1.1.} Identify a point $\vecu \in [0, 1]^d$ with the indicator function $\I_{(-\vc{\infty}, \vecu]}$ on $\mathbb{R}^d$. The resulting class $\mathcal{F} = \{ \I_{(-\vc{\infty}, \vecu]} : \vecu \in [0, 1]^d \}$, being bounded (by $1$) and VC \cite[Example~2.6.1]{Vaart-Wellner}, it satisfies the uniform entropy condition~(2.5.1) in \cite{Vaart-Wellner}; see Theorem~2.6.7 in the same book.  From their Theorem~2.8.3, we obtain the uniform Donsker property
\begin{equation}
\label{eq:uniformDonsker}
\sup_{P \in \mathcal{P}} \sup_{h \in \operatorname{BL}_1}
\left\lvert \E_P^*[ h( \bbU_{n,P} ) ] - \E[ h( \bbU^P ) ] \right\rvert
\tendsto 0, \qquad n \to \infty.
\end{equation}
The supremum over $h$ in \eqref{eq:1.1.1} is bounded by the double supremum over $P$ and $h$ in \eqref{eq:uniformDonsker}. The convergence in \eqref{eq:1.1.1} is thus proved.
\smallskip


\emph{Step~1.1.2.} 
We need to show that, almost surely, $\bbU^{P_n^\beta} \rightsquigarrow \bbU^C$ as $n \to \infty$. All processes involved are tight, centered Gaussian processes, with covariance functions determined in \eqref{eq:UC} using $\bbC_n^P$ or $C$. The strong consistency of the empirical copula, together with \cite[Proposition~2.8]{SST2017}, yields 
\begin{equation} 
\label{eq:strongconsistency}
\lVert \bbC_n^\beta - C \rVert_\infty \to 0, 
\qquad \text{$n \to \infty$, a.s.} 
\end{equation}
This property implies \eqref{eq:1.1.2}. First, \eqref{eq:strongconsistency} implies the almost sure convergence of the covariance function of $\bbU^{P_n^\beta}$ to that of $\bbU^C$, and thus the almost sure convergence of the finite-dimensional distributions. Second, the asymptotic tightness a.s.\ follows from the uniform continuity of the trajectories with respect to their respective intrinsic standard deviation semimetrics \eqref{eq:sdmetric} and the uniform convergence a.s.\ of these standard deviation semimetrics, again by \eqref{eq:strongconsistency}.
\smallskip

\emph{Step 1.2} The asymptotic measurability of $\gamma_n^{\#}$ follows from the \emph{unconditional} (i.e., jointly in $\vecX_1,\ldots,\vecX_n,\vecV_1^{\#},\ldots,\vecV_n^{\#}$) weak convergence $\gamma_n^{\#} \rightsquigarrow \bbU^C$ as $n \to \infty$. This claims can be divided into the convergence of the finite-dimensional distributions and asymptotic tightness. The former can be shown using the Lindeberg central limit theorem for triangular arrays conditional on $\vecX_1, \ldots, \vecX_n$, using a similar method to that in the proof of Theorem~23.4 in~\cite{Vaart1998}. The latter follows as in Theorems~2.5.2 and~2.8.3 in~\cite[p.~128 and~171]{Vaart-Wellner}, conditional on $\vc{X}_1, \ldots, \vc{X}_n$ using the fact that the class of indicator functions of cells in $\bbR^d$ is a VC-class \cite[Example~2.5.4]{Vaart-Wellner}.
\smallskip

\emph{Step 2.}
Consider a map $\Phi$ that sends a cumulative distribution function $H$ on $[0, 1]^d$ whose marginals do not assign mass at zero to the function $\vecu \mapsto H(H_1^-(u_1),\ldots,H_d^-(u_d))$. We have $\cop^{\#} = \Phi(\bbG_n^{\#})$ and $\copb = \Phi(\copb)$. By \cite[Theorem~2.4]{Buech-Volg2013}, the map $\Phi$ is Hadamard differentiable at the true copula $C$ tangentially to a certain set $\mathbb{D}_0$ at which the distribution of $\mathbb{U}^C$ is concentrated.  By \eqref{eq:toshow:gamma}, the form of the Hadamard derivative $\Phi_C'$ of $\Phi$ at $C$ together with the functional delta method for the bootstrap \cite[Theorem~3.9.11]{Vaart-Wellner} yield conditional weak convergence in probability of $\alpha_n^{\#} = \sqrt{n}\{\Phi(\bbG_n^{\#}) - \Phi(\copb)\}$ to $\bbG^C = \Phi_C'(\bbU^C)$.

Since $\lvert \tilbbC_n - \cop \rvert \le d/n$ and $\lvert \tilbbC_n^{\#} - \cop^{\#} \rvert \leq d/n$ by \eqref{ineq:emp-cop:rank-cop}, we obtain the conditional weak convergence in probability of $\tilde{\alpha}_n^{\#}$ to $\bbG^C$. 

Finally, since $\alpha_n^{\beta \#}(\vc{u}) = \E_S[ \tilde{\alpha}_n^{\#}(S_1/n, \ldots, S_d/n) ]$ as in \eqref{eq:alphaES}, we arrive at the conditional weak convergence in probability of $\alpha_n^{\beta \#}$ to $\bbG^C$ in a way similar to the proof of Proposition~\ref{prop:asymp-equiv}.
\qed

\bibliographystyle{plain}
\bibliography{tsuka_ref}

\begin{thebibliography}{10}

\bibitem{Buech-Dette2010}
A.~B{\"{u}}cher and H.~Dette.
\newblock A note on bootstrap approximations for the empirical copula process.
\newblock {\em Statistics and Probability Letters}, 80:1925--1932, 2010.

\bibitem{Buech-Volg2013}
A.~B{\"{u}}cher and S.~Volgushev.
\newblock Empirical and sequential empirical copula processes under serial
  dependence.
\newblock {\em Journal of Multivariate Analysis}, 119:61--70, 2013.

\bibitem{Deheu79}
P.~Deheuvels.
\newblock La fonction de d\'ependence empirique et ses propri\'et\'es, un test
  non param\'etrique d'ind\'ependance.
\newblock {\em Bulletin de la classe des sciences, Acad\'emie Royale de
  Belgique, 5e s\'erie}, 65:274--292, 1979.

\bibitem{Efron1982siam}
B.~Efron.
\newblock {\em The Jackknife, the Bootstrap and Other Resampling Plans}.
\newblock Society for Industrial and Applied Mathematics, Philadelphia, 1982.

\bibitem{Fer-Rad-Weg04}
J.-D. Fermanian, D.~Radulovi\'{c}, and M.~J. Wegkamp.
\newblock Weak convergence of empirical copula processes.
\newblock {\em Bernoulli}, 10:847--860, 2004.

\bibitem{Gen-Gho-Riv95}
C.~Genest, K.~Ghoudi, and L.-P. Rivest.
\newblock A semiparametric estimation procedure of dependence parameters in
  multivariate families of distributions.
\newblock {\em Biometrika}, 82:543--552, 1995.

\bibitem{Gen-Nes-Que2012}
C.~Genest, J.~Ne\v{s}lehov\'{a}, and J.-F. Quessy.
\newblock Tests of symmetry for bivariate copulas.
\newblock {\em Annals of the Institute of Statistical Mathematics},
  64:811^^e2^^80^^93--834, 2012.

\bibitem{Hall92}
P.~Hall.
\newblock {\em The Bootstrap and Edgeworth Expansion}.
\newblock Springer-Verlag, New York, 1992.

\bibitem{copulaR}
M.~Hofert, I.~Kojadinovic, M.~Maechler, J.~Yan, and J.~G. Ne\v{s}lehov\'{a}.
\newblock Package `copula', \textsf{R} package version 0.999-19, 2018.

\bibitem{Hol-Wol-Chi2014}
M.~Hollander, D.~A. Wolfe, and E.~Chicken.
\newblock {\em Nonparametric Statistical Methods}.
\newblock John Wiley \& Sons, Hoboken, New Jersey, third edition, 2014.

\bibitem{Khoudraji1995}
A.~Khoudraji.
\newblock {\em Contributions \`{a} l'\'{e}ude des copules et \`{a} la
  mod\'{e}lisation de valeurs extr\^{e}mes bivari\'{e}s}.
\newblock PhD thesis, Universit\'{e} Laval, Qu\'{e}bec, Canada, 1995.

\bibitem{Kiriliouk-Segers-Tafakori2018}
A.~Kiriliouk, J.~Segers, and L.~Tafakori.
\newblock An estimator of the stable tail dependence functionbased on the
  empirical beta copula.
\newblock {\em Extremes}, 21:581--600, 2018.

\bibitem{Kojadinovic-Yan2012}
I.~Kojadinovic and J.~Yan.
\newblock A non-parametric test of exchangeability for extreme-value and
  left-tail decreasing bivariate copulas.
\newblock {\em Scandinavian Journal of Statistics}, 39:480--496, 2012.

\bibitem{Kosorok2008}
M.~R. Kosorok.
\newblock {\em Introduction to Empirical Processes and Semiparametric
  Inference}.
\newblock Springer, New York, 2008.

\bibitem{Nelsen2006}
R.~B. Nelsen.
\newblock {\em An Introduction to Copulas}.
\newblock Springer-Verlag, New York, second edition, 2006.

\bibitem{Remillard-Scaillet2009}
B.~R\'{e}millard and O.~Scaillet.
\newblock Testing for equality between two copulas.
\newblock {\em Journal of Multivariate Analysis}, 100:377--386, 2009.

\bibitem{Segers2012}
J.~Segers.
\newblock Asymptotics of empirical copula processes under nonrestrictive
  smoothness assumptions.
\newblock {\em Bernoulli}, 18:764--782, 2012.

\bibitem{SST2017}
J.~Segers, M.~Sibuya, and H.~Tsukahara.
\newblock The empirical beta copula.
\newblock {\em Journal of Multivariate Analysis}, 155:35--51, 2017.

\bibitem{Shao-Tu95}
J.~Shao and D.~Tu.
\newblock {\em The Jackknife and Bootstrap}.
\newblock Springer-Verlag, New York, 1995.

\bibitem{Sklar59}
M.~Sklar.
\newblock Fonctions de r\'epartition \'a n dimensions et leurs marges.
\newblock {\em Publ.\ Inst.\ Statist.\ Univ.\ Paris}, 8:229--231, 1959.

\bibitem{Tsuka05}
H.~Tsukahara.
\newblock Semiparametric estimation in copula models.
\newblock {\em Canadian Journal of Statistics}, 33:357--375, 2005.
\newblock [Erratum: \textit{Canadian Journal of Statistics}, \textbf{39},
  734--735 (2011)].

\bibitem{Vaart1998}
A.~W. van~der Vaart.
\newblock {\em Asymptotic Statistics}.
\newblock Cambridge University Press, Cambridge, 1998.

\bibitem{Vaart-Wellner}
A.~W. van~der Vaart and J.~A. Wellner.
\newblock {\em Weak Convergence and Empirical Processes: With Applications to
  Statistics}.
\newblock Springer-Verlag, New York, 1996.

\end{thebibliography}

\end{document}